\documentclass[submission]{FPSAC2018}

\newtheorem{defi}{Definition}
\newtheorem{thm}{Theorem}
\newtheorem{prop}{Proposition}
\newtheorem{lem}{Lemma}
\newtheorem{cor}{Corollary}

\title[Cylindric RPP and 2D TQFT]{Cylindric Reverse Plane Partitions and 2D TQFT}

\author{Christian Korff\thanks{\href{mailto:christian.korff@glasgow.ac.uk}{christian.korff@glasgow.ac.uk}}\addressmark{1}  \and David Palazzo\thanks{\href{mailto:d.palazzo1@research.gla.ac.uk}{d.palazzo1@research.gla.ac.uk}. David Palazzo is supported by a College of Science and Engineering Scholarship of the University of Glasgow.}\addressmark{1}}

\address{\addressmark{1}School of Mathematics and Statistics, University of Glasgow, Glasgow G12 8QQ}

\received{14 November 2017}

\abstract{The ring of symmetric functions carries the structure of a Hopf algebra. When computing the coproduct of complete symmetric functions $h_\lambda$ one arrives at weighted sums over reverse plane partitions (RPP) involving binomial coefficients. Employing the action of the extended affine symmetric group at fixed level $n$ we generalise these weighted sums to cylindric RPP and define cylindric complete symmetric functions. The latter are shown to be $h$-positive, that is, their expansions coefficients in the basis of complete symmetric functions are non-negative integers. We state an explicit formula in terms of tensor multiplicities for irreducible representations of the generalised symmetric group. Moreover, we relate the cylindric complete symmetric functions to a 2D topological quantum field theory (TQFT) that is a generalisation of the celebrated $\mathfrak{\widehat{sl}}_n$-Verlinde algebra or Wess-Zumino-Witten fusion ring, which plays a prominent role in the context of vertex operator algebras and algebraic geometry.}

\keywords{reverse plane partitions, symmetric functions, topological quantum field theory, Verlinde algebra, generalised symmetric group }

\begin{document}

\maketitle

\section{Symmetric functions and reverse plane partitions}
This section recalls some known results about the ring of symmetric functions and reverse plane partitions. It serves as a motivation for our definition of weighted sums over cylindric reverse plane partitions in Section 2 and the definition of a certain class of 2D TQFTs in Section 3, which are the new results. Detailed proofs and a full discussion will be published elsewhere \cite{KorffPalazzo}, here we simply present a summary of our results.

\subsection{Pieri type rules and coproduct}

Let $\Lambda=\mathbb{Z}[h_1,h_2,\ldots]$ denote the ring of symmetric functions equipped with the Hall inner product $\langle m_\lambda,h_\mu\rangle =\delta_{\lambda\mu}$, where $m_\lambda$ denotes the monomial symmetric function and $h_\mu=h_{\mu_1}h_{\mu_2}\cdots$ the complete symmetric function with $h_r=\sum_{1\leq i_1\leq\cdots\leq i_r}x_{i_1}x_{i_2}\cdots x_{i_r}$ and $\lambda,\mu\in\mathcal{P}$ partitions.  In particular, each of the sets $\{m_\lambda\}_{\lambda\in\mathcal{P}}$ and $\{h_\lambda\}_{\lambda\in\mathcal{P}}$ forms a $\mathbb{Z}$-basis of the ring $\Lambda$.

Introduce a coproduct $\Delta: \Lambda \to \Lambda \otimes \Lambda$ via the relation $\langle \Delta f, g \otimes h \rangle = \langle f, gh \rangle$ for $f,g,h\in\Lambda$; compare with \cite[Ch.I.5, Ex. 25]{macdonald1998symmetric}. The latter allows one to use the product formulae
\begin{equation} \label{pierirule} 
m_\lambda m_\mu=\sum_{\nu}f_{\lambda\mu}^{\nu}m_{\nu}
\qquad\text{and}\qquad h_\lambda m_\mu=\sum_{\nu}\chi_{\nu/\mu}(\lambda)m_\nu
\end{equation}
to define the `skew complete symmetric function' $h_{\lambda/\mu}$ via
\begin{equation} \label{hlambdamu}
\Delta h_\lambda=\sum_{\mu}h_{\lambda/\mu}\otimes h_\mu,\qquad
 h_{\lambda/\mu}=\sum_{\nu}f_{\mu\nu}^\lambda h_\nu
 =\sum_{\nu}\chi_{\lambda/\mu}(\nu)m_\nu
\;.
\end{equation}
The expansion coefficients in both formulae are non-negative integers and have combinatorial expressions, which we know recall. 

For the product of monomial symmetric functions the expansion coefficients $f_{\lambda\mu}^\nu$ in \eqref{pierirule} are given by the cardinality of the set $\{(\alpha,\beta)~|~\alpha+\beta=\nu\}$ where respectively $\alpha$ and $\beta$ are distinct permutations of $\lambda$ and $\mu$; see e.g. \cite{butler1993nonnegative}. This also fixes the expansion coefficients for the second product in \eqref{pierirule}. Recall the known identity $h_\lambda=\sum_{\mu}L_{\lambda\mu}m_\mu$, where $L_{\lambda\mu}$ is the number of $\mathbb{N}$-matrices $A=(a_{ij})_{i,j\geq1}$ which have only a finite number of nonzero entries and row sums $\text{row}(A)=\lambda$ and column sums $\text{col}(A)=\mu$; see e.g. \cite[Ch.7.5, Prop 7.5.1]{stanley_fomin_1999}. Then it follows at once that
\begin{equation}\label{skewChi}
\chi_{\nu/\mu}(\lambda)=\sum_{\alpha}L_{\lambda\alpha}f^{\nu}_{\alpha\mu}\;.
\end{equation}
Here we are interested in an alternative expression using reverse plane partitions as the latter suggests a natural generalisation to cylindric reverse plane partitions.

\subsection{Weighted sums over reverse plane partitions}

Given two partitions $\lambda,\mu$ with $\mu\subset \lambda$ recall that a {\em reverse plane partition} (RPP) $\pi$ of skew shape $\lambda/\mu$ is a sequence $\{\lambda^{(i)}\}_{i=0}^l$ of partitions $\lambda^{(i)}$ with %
$
\mu=\lambda^{(0)}\subset\lambda^{(1)}\subset\cdots\subset\lambda^{(l)}=\lambda 
$. %
As usual we refer to the vector $\theta(\pi)=(|\lambda^{(1)}/\lambda^{(0)}|,|\lambda^{(2)}/\lambda^{(1)}|,\ldots)$ as the {\em weight} of $\pi$ and denote by $x^\pi$ the monomial $x^\pi=x_1^{\theta_1(\pi)}x_2^{\theta_2(\pi)}\cdots$ in the indeterminates $x_i$. 
Alternatively, we can think of a RPP as a map $\pi:\lambda/\mu\to\mathbb{N}$ which assigns to the squares $s=(x,y)\in \lambda^{(i)}/\lambda^{(i-1)}\subset \mathbb{Z}\times\mathbb{Z}$ the integer $i$. The result is a skew tableau whose entries are non-decreasing along each row from left to right and down each column. 

We now introduce weighted sums over RPP in terms of binomial coefficients. Given any pair of partitions $\lambda,\mu\in\mathcal{P}$ denote by $\chi_{\lambda/\mu}$ the number of distinct permutations $\alpha$ of $\mu$ such that $\alpha\subset\lambda$. The latter has the following explicit expression in terms of the conjugate partitions $\lambda',\mu'$,
\begin{equation}\label{chi} 
 \chi_{\lambda / \mu} = \begin{cases}
                          \prod_{i \ge 1} {{\lambda'_i - \mu'_{i+1}} \choose {\mu'_i - \mu'_{i+1}}},
                          & \mu \subset \lambda \\
                          0, & \text{otherwise}
                         \end{cases}\;.
\end{equation}
The following result is probably known to experts but we were unable to find it in the literature.
\begin{lem}
The skew complete symmetric function defined in \eqref{hlambdamu} is the weighted sum
 \begin{equation}
  h_{\lambda / \mu}(x) = \sum_{\pi} \chi_\pi \, x^\pi,
  \qquad \chi_\pi=\prod_{i\ge 1} \chi_{\lambda^{(i)}/\lambda^{(i-1)}}
 \end{equation}
over all reverse plane partitions $\pi$ of shape $\lambda / \mu$. In particular, the coefficient \eqref{skewChi} has the alternative expression $\chi_{\lambda/\mu}(\nu)=\sum_{\pi} \chi_\pi$, where the sum is over all reverse plane partitions $\pi$ of shape $\lambda/\mu$ and weight $\theta(\pi)=\nu$.
\end{lem}

\subsection{The expansion coefficients in terms of the symmetric group}

We now project onto the ring $\Lambda_k=\mathbb{Z}[x_1,\ldots,x_k]^{S_k}$ of symmetric functions in $k$-variables by setting $x_{k+1}=x_{k+2}=\cdots =0$. Here $S_k$ is the symmetric group in $k$ letters. %
Let $\mathcal{P}_k = \bigoplus_{i=1}^k \mathbb{Z} \epsilon_i$ be the $\mathfrak{gl}_k$
weight lattice with standard basis $\epsilon_1,\dots,\epsilon_k$ and inner product $(\epsilon_i,\epsilon_j)=\delta_{ij}$. Denote by $\mathcal{P}_k^+\subset\mathcal{P}_k$ the {\em positive dominant weights}. 

Each $\lambda\in\mathcal{P}_k$ defines a map $\lambda:[k]\to\mathbb{Z}$ in the obvious manner and we shall consider the right action $\mathcal{P}_k\times S_k\to\mathcal{P}_k$ given by $(\lambda,w)\mapsto \lambda\circ w=(\lambda_{w(1)},\ldots,\lambda_{w(k)})$. For a fixed weight $\mu$ denote by $S_\mu\subset S_k$ its stabiliser group. The latter has cardinality $|S_\mu|=\prod_{i\in\mathbb{Z}}m_{i}(\mu)!$ with $m_i(\mu)$ being the number of elements in $[k]$ which are mapped to $i$ under $\mu$. Given any permutation $w\in S_k$ there exists a unique decomposition $w=w_\mu w^\mu$ with $w_\mu\in S_\mu$ and $w^\mu$ a minimal length representative of the right coset $S_\mu w$. Denote by $S^\mu\subset S_k$ the set of all minimal length coset representatives in $S_\mu\backslash S_k$.
\begin{lem}
Let $\lambda,\mu,\nu\in\mathcal{P}_k^+$. (i) The expansion coefficient $f_{\lambda\mu}^\nu$ in \eqref{pierirule} can be expressed as the cardinality of the set
\begin{equation}\label{setf}
\{(w,w')\in S^\lambda\times S^\mu~|~\lambda\circ w+\mu\circ w' =\nu\}\;.
\end{equation}
(ii) The specialisation of the weight factor $\chi_{\lambda/\mu}$ defined in \eqref{chi} to elements in $\mathcal{P}_k^+$ equals the cardinality of the set
 \begin{equation} \label{setChi}
   \{ w \in S^\mu~|~ \mu\circ w \le \lambda \}\,,
 \end{equation}
where $\mu\circ w\le \lambda$ is shorthand notation for $\mu_{w(i)} \le \lambda_i$ for all $i\in [k]$. 
\end{lem}

\section{Cylindric Reverse Plane Partitions}

Cylindric plane partitions were first considered by Gessel and Krattenthaler in \cite{gessel1997cylindric}. There has been a growing interest in cylindric Schur functions \cite{mcnamara2006cylindric} and their connections with other areas in mathematics. %
%
Here we introduce {\em cylindric complete symmetric functions} as weighted sums over cylindric RPP by generalising the definition of the sets \eqref{setf} and \eqref{setChi} to the extended affine symmetric group $\hat S_k=\mathcal{P}_k\rtimes S_k$. In the last section we will relate the expansions of these cylindric complete symmetric functions in the $h_\lambda$-basis to the (non-negative) structure constants of generalised Verlinde algebras. 

\subsection{Lusztig's realisation of the affine symmetric group}
Recall the realisation of the affine symmetric group $\tilde{S}_k=\mathcal{Q}_k\rtimes S_k$, where $\mathcal{Q}_k\subset\mathcal{P}_k$ is the root lattice, in terms of bijections $\tilde w:\mathbb{Z}\to\mathbb{Z}$ that are subject to the following two conditions: 
\begin{equation}\label{Lusztig}
\tilde w(m+k)=\tilde w(m)+k,\;\forall m\in\mathbb{Z}\qquad\text{and}\qquad\sum_{i=1}^k {\tilde w(m)} =\binom{k}{2}\;.
\end{equation}
As in the finite case the group multiplication is given by composition. This realisation of $\tilde S_k$ first appeared in \cite{lusztig1983some} and has subsequently used by other authors \cite{bjorner1996affine}. The group is generated by the simple Weyl reflections $\{ \sigma_0, \sigma_1, \ldots,\sigma_{k-1} \}$ which as maps $\mathbb{Z}\to\mathbb{Z}$ are defined as
\begin{equation}
\sigma_i(m) = 
\begin{cases}
m+1, & m=i \text{ mod } k \\
m-1, & m=i+1 \text{ mod } k \\
m, & \text{otherwise}
\end{cases}
\end{equation}
and one verifies that they satisfy the familiar relations
 \begin{eqnarray*}
  \sigma_i\circ\sigma_i=\operatorname{Id}, \qquad  \sigma_i\circ \sigma_{i+1}\circ \sigma_i = \sigma_{i+1}\circ \sigma_i\circ \sigma_{i+1}, \qquad
  \sigma_i\circ \sigma_j = \sigma_j\circ \sigma_i \quad\text{ for } |i-j|>1\,,
 \end{eqnarray*}
where all indices are understood modulo $k$. 

\subsection{The extended affine symmetric group}

We now state a realisation of the {\em extended} affine symmetric group in terms of bijections $\mathbb{Z}\to\mathbb{Z}$. Define the shift operator $\tau:\mathbb{Z}\to\mathbb{Z}$ by $m\mapsto\tau(m)=m-1$. Then one has the identities 
\begin{equation}
\tau\circ\sigma_{i+1}=\sigma_i\circ\tau,\qquad i=0,1,\ldots,k-1\;.
\end{equation}
That is, the group generated by $\langle\tau,\sigma_0, \sigma_1, \ldots,\sigma_{k-1} \rangle$ is the extended affine symmetric group $\hat S_k=\mathcal{P}_k\rtimes S_k$. To see this more clearly, we introduce the additional generators $y_k=\tau\circ\sigma_1\circ\sigma_2\circ\cdots\circ\sigma_{k-1}$ and $y_i=\sigma_i\circ y_{i+1}\circ\sigma_i$ for $i=1,2,\ldots,k-1$. Then any element $\hat w\in\hat S_k$ can be written as $\hat w=w\circ y^\lambda$, where $\lambda\in\mathcal{P}_k$, $y^\lambda=y_1^{\lambda_1}\circ\cdots\circ y_k^{\lambda_k}$ and $w\in S_k$. Note that $\hat w:\mathbb{Z}\to\mathbb{Z}$ still satisfies the first property in \eqref{Lusztig} but not the second. Because any element $\hat w\in\hat S_k$ can be expressed as $\hat w= \tilde w\circ\tau^d$ for some $d\in\mathbb{Z}$ and $\tilde w\in\tilde S_k$, we have the following realisation of $\hat S_k$ generalising \eqref{Lusztig}:
\begin{lem}
The extended affine symmetric group $\hat S_k$ can be realised as the bijections $\hat w:\mathbb{Z}\to\mathbb{Z}$ subject to the conditions
\begin{equation}\label{extS}
\hat w(m+k)=\hat w(m)+k,\;\forall m\in\mathbb{Z}\qquad\text{and}\qquad\sum_{m=1}^k {\hat w(m)} =\binom{k}{2}\mod k\;.
\end{equation}
\end{lem} 
Our main interest in this realisation of $\hat S_k$ is that it naturally leads to the consideration of cylindric loops.

\subsection{Cylindric loops and cylindric reverse plane partitions}

Fix $n\in\mathbb{N}$. We are now generalising the notion of the weight lattice in order to define a level-$n$ action of the extended affine symmetric group. Let $\mathcal{P}_{k,n}$ denote the set of functions $\lambda:\mathbb{Z}\to\mathbb{Z}$ subject to the constraint $\lambda_{i+k}=\lambda_i-n$ for all $i\in\mathbb{Z}$. Then the map $\mathcal{P}_{k,n}\times\hat S_k\to \mathcal{P}_{k,n}$ with $(\lambda,\hat w)\mapsto \lambda\circ\hat w$ defines a right action.  One can convince oneself that this is the familiar level-$n$ action of $\hat S_k$ on the weight lattice $\mathcal{P}_k$ by observing that each $\lambda\in\mathcal{P}_{k,n}$ is completely fixed by its values $(\lambda_1,\ldots,\lambda_k)$ on the set $[k]$. In particular, the `alcove'
\begin{equation}\label{alcove0}
 \mathcal{A}_{k,n} = \{ \lambda\in\mathcal{P}_{k,n} ~|~n\ge \lambda_1 \ge \lambda_2 \ge
 \dots \ge \lambda_k > 0 \}
\end{equation}
 is a fundamental domain with respect to this level-$n$ action of $\hat S_k$ on $\mathcal{P}_{k,n}$. That is, for any $\lambda\in\mathcal{P}_{k,n}$ the orbit $\lambda\hat S_k$ intersects $\mathcal{A}_{k,n}$ in a unique point.
 
Given $\lambda\in\mathcal{A}_{k,n}$ and $d\in\mathbb{Z}$ denote by $\lambda[d]$ the (doubly) infinite sequence
 \[
 \lambda[d]=( \dots,\lambda[d]_{-1},\lambda[d]_0,\lambda[d]_1,\dots)=
 ( \dots,\lambda_{-d-1},\lambda_{-d},\lambda_{1-d},\dots)\,,
 \]
that is, the image $\lambda\circ\tau^d(\mathbb{Z})$ of the map $\lambda\circ\tau^d:\mathbb{Z}\to\mathbb{Z}$. This sequence defines a lattice path in $\mathbb{Z}\times\mathbb{Z}$ which repeats itself and is therefore called a {\em cylindric loop}. A {\em cylindric skew diagram} or {\em cylindric shape} can now be defined as the number of lattice points between two cylindric loops: let  $\lambda,\mu \in \mathcal{A}_{n,k}$ be such that $\mu_i \le \lambda_{i-d}=(\lambda\circ\tau^d)_i$ for all $i \in \mathbb{Z}$, then we write
$\mu[0] \le \lambda[d]$ and say that the set 
 \[
 \lambda / d / \mu = \{ (i,j) \in \mathbb{Z}^2 ~|~ \mu[0]_i\le j \le \lambda[d]_i  \}
 \]
 is a cylindric skew diagram of degree $d$.
\begin{defi}
A {\em cylindric reverse plane partition} of shape $\Theta=\lambda/d/\mu$ is a map $\hat\pi:\Theta \to \mathbb{N}$ such that for any $(i,j) \in \Theta$ one has  $\hat\pi(i,j) = \hat\pi(i+k,j-n)$ together with
 \[
  \hat\pi(i,j) \le \hat\pi(i+1,j)\quad\text{ and }\quad\hat\pi(i,j) \le \hat\pi(i,j+1),\quad \text{ if } (i+1,j),(i,j+1) \in \Theta \;.
\]
In other words, the entries in the squares between the cylindric loops $\mu[0]$ and $\lambda[d]$ are non-decreasing from left to right in rows and down columns.
\end{defi}

\begin{figure}
\centering
\includegraphics[width=0.4\textwidth]{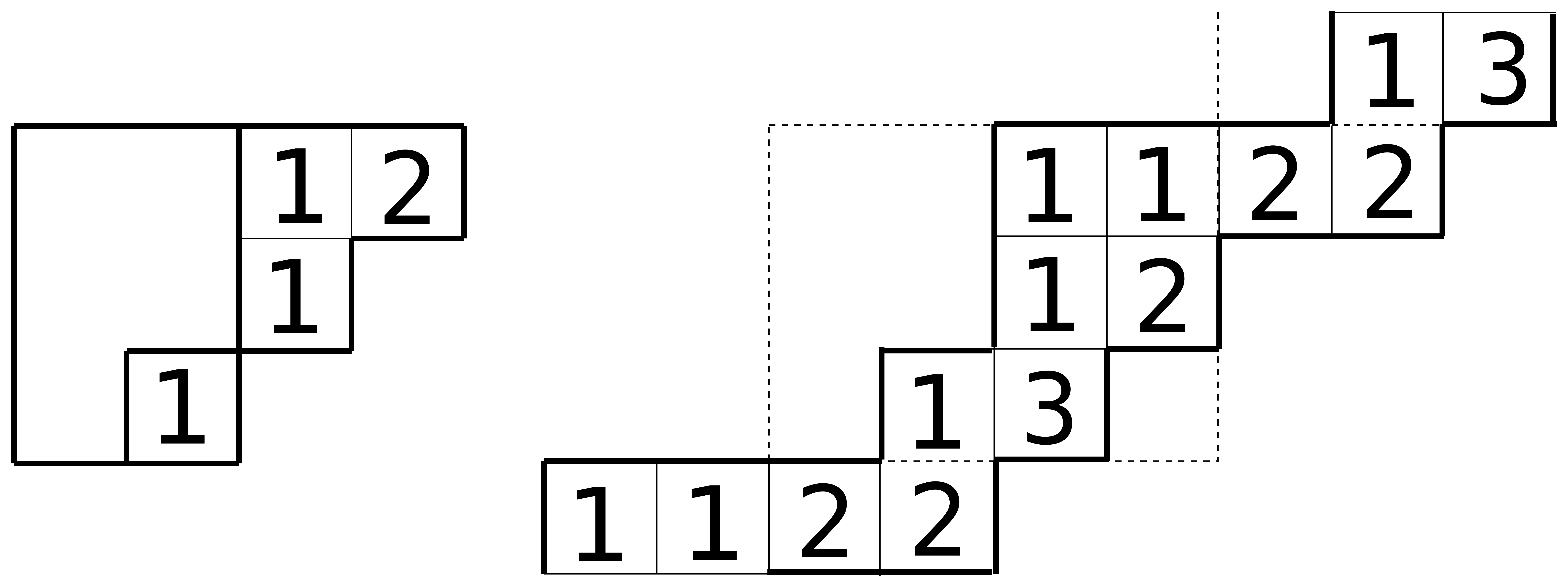}
\caption{On the left a RPP of shape $(4,3,2)/(2,2,1)$ and weight $\theta=(3,1,0)$ and on the right a cylindric RPP of shape $(4,3,2)/1/(2,2,1)$ and weight $\theta=(4,3,1)$ with $n=4$ and $k=3$.}
\label{fig:plot}
\end{figure}

Alternatively, $\hat\pi$ can be defined as a sequence of cylindric loops 
\begin{equation} \label{SequenceCylLoops}
( \lambda^{(0)}[0]=\mu[0], \lambda^{(1)}[d_1], \dots, \lambda^{(l)}[d_l]=\lambda[d] )
\end{equation}
with $\lambda^{(i)} \in \mathcal{A}_{k,n}$ and $ d_i - d_{i-1} \ge 0$ such that 
$\hat\pi^{-1}(i)=\lambda^{(i)} / (d_i-d_{i-1}) / \lambda^{(i-1)}$ is a cylindric skew diagram; see Figure \ref{fig:plot} for a simple example. %
The weight of $\hat\pi$ is the vector $\theta(\hat\pi)=(\theta_1,\ldots,\theta_l)$ where $\theta_i$ is the number of lattice points $(a,b)\in \lambda^{(i)} / (d_i-d_{i-1}) / \lambda^{(i-1)}$ with $1\leq a\leq k$. 

\subsection{Cylindric complete symmetric functions}

Given $\mu\in\mathcal{A}_{k,n}$ note that $\mu_1-\mu_k<n$ and, hence, its stabiliser group $S_\mu\subset S_k\subset\hat S_k$. Define $\tilde S^{\mu}$ as the minimal length representatives of the cosets $S_\mu\backslash\tilde S_k$. Consider now the conjugate partition $\mu'$, which is defined by $\mu'_i-\mu'_{i+1}=m_i(\mu)$ for $i=1,\ldots,n$, where $m_i(\mu)$ is the multiplicity of the part $i$. Then define $\mu':\mathbb{Z}\to\mathbb{Z}$ by setting $\mu'_{i+n}=\mu'_i-k$. 
\begin{lem}
For $\lambda,\mu \in \mathcal{A}_{k,n}$ and $d \in \mathbb{Z}_{\ge 0}$, let $\chi_{\lambda/d/\mu}$ be the cardinality of the set
\begin{equation} \label{SetAffinePermChi}
  \{ \tilde{w}\in \tilde{S}^\mu~|~\mu\circ\tilde{w} \le\lambda\circ\tau^d \}\;.
\end{equation}
Then we have the following expression in terms of binomial coefficients and cylindric loops,
\begin{equation} \label{CylindricChi}
 \chi_{\lambda/d/\mu}=\prod_{j=1}^n { {(\lambda\circ\tau^d)'_j-\mu'_{j+1}} \choose {\mu'_j-\mu'_{j+1}} } - 
  \prod_{j=1}^n { {(\lambda\circ\tau^{d-1})'_j-\mu'_{j+1}} \choose {\mu'_j-\mu'_{j+1}} }\;.
\end{equation}
Here the binomial coefficients are defined to be zero if any of their arguments is negative.
\end{lem}
As in the non-cylindric case \eqref{chi} we employ \eqref{CylindricChi} to define weighted sums over cylindric RPP: given $\hat\pi$ set $\chi_{\hat \pi}=\prod_{i\ge 1} \chi_{\lambda^{(i)}/(d_i-d_{i-1})/\lambda^{(i-1)}}$, where the cylindric skew diagram $\lambda^{(i)}/(d_i-d_{i-1})/\lambda^{(i-1)}$ is the pre-image $\hat\pi^{-1}(i)$. 

\begin{defi}
For $\lambda,\mu \in \mathcal{A}_{k,n}$ and $d \in \mathbb{Z}_{\ge 0}$, introduce the {\em cylindric complete symmetric function} $h_{\lambda/d/\mu}$ as the weighted sum
\begin{equation}\label{cylh}
h_{\lambda/d/\mu}(x)=\sum_{\hat\pi} \chi_{\hat\pi} x^{\hat\pi}
\end{equation}
over all cylindric reverse plane partitions $\hat\pi$ of shape $\lambda/d/\mu$.
\end{defi}
Note that the cylindric complete symmetric functions introduced here are different from the cylindric skew Schur functions considered in \cite{mcnamara2006cylindric} through the inclusion of the factor $\chi_{\hat\pi}$ and involve sums over binomial coefficients.
\begin{thm}
Let $\lambda,\mu\in\mathcal{A}_{k,n}$ and $d\in\mathbb{Z}_{\ge 0}$. The function $h_{\lambda/d/\mu}$ is symmetric and one has the expansion
\begin{equation}\label{cylh2h}
h_{\lambda/d/\mu}=\sum_{\nu\in\mathcal{P}^+_{k}}N_{\mu\nu}^{\lambda}
\;h_\nu
\end{equation}
into the basis $\{h_\nu\}_{\nu\in\mathcal{P}}$ of complete symmetric functions in $\Lambda$, where 
\begin{itemize}
\item[(i)] the sum is restricted to those $\nu\in\mathcal{P}^+_k$ for which $|\nu|=nd+|\lambda|-|\mu|$ and
\item[(ii)]  $N_{\mu\nu}^{\lambda}\in\mathbb{Z}_{\geq 0}$ are given by the cardinality of the set
\begin{equation}\label{setN}
\{(w,w')\in S^\mu\times S^\nu~|~\mu\circ w+\nu\circ w'=\lambda\circ y^\alpha\text{ for some }\alpha\in\mathcal{P}_k\text{ with }\;|\alpha|=d\},
\end{equation}
with $y_i\in\hat S_k$ being the translations in the weight lattice introduced earlier.
\end{itemize}
\end{thm}
Note that when setting $d=0$ we recover the (non-cylindric) skew complete cylindric function, $h_{\lambda/0/\mu}=h_{\lambda/\mu}$. In particular, $N_{\lambda\mu}^\nu=f_{\lambda\mu}^\nu$ as defined in \eqref{setf}, since it follows from $\lambda,\mu\in\mathcal{A}_{k,n}$ that $\alpha_i=0$ for all $i=1,\ldots,k$.
\begin{lem}
Let $\lambda,\mu\in\mathcal{A}_{k,n}$ and $\nu\in\mathcal{P}^+_k$. Denote by $\check\nu\in\mathcal{A}_{k,n}$ the unique intersection point of the orbit $\nu\hat S_k$ with the alcove \eqref{alcove0}. Then the expansion coefficients for $\nu$ in \eqref{cylh2h} can be re-written in terms of those for $\check\nu$ and multinomial coefficients,
\begin{equation}\label{reducedN}
N_{\mu\nu}^{\lambda}=N_{\mu\check\nu}^{\lambda}
\binom{m_n(\check\nu)}{m_0(\nu),m_{n}(\nu),m_{2n}(\nu),\ldots}
\prod_{i=1}^{n-1}\binom{m_i(\check\nu)}{m_i(\nu),m_{i+n}(\nu),m_{i+2n}(\nu),\ldots}\;,
\end{equation}
where $m_j(\nu)$ and $m_j(\check\nu)$ are the multiplicities of the part $j$ in $\nu$ and $\check\nu$, respectively.
\end{lem}

\section{The generalised symmetric group and 2D TQFT}

In view of the product and coproduct formulae \eqref{pierirule} in $\Lambda$, one might ask whether similar formulae hold for the cylindric complete symmetric functions \eqref{cylh} in an appropriate ring or algebra. We now construct such an algebra as a particular quotient of $\Lambda_k$ such that the expansion coefficients \eqref{cylh2h} are the structure constants with respect to the (projected) monomial symmetric functions. The quotient is finite-dimensional and can be endowed with the structure of a Frobenius algebra (more precisely a Frobenius extension).  Based on work of Atiyah \cite{atiyah1988topological} Frobenius algebras are categorically equivalent to 2D topological quantum field theories which are functors from the category of 2-cobordisms to finite-dimensional vector spaces; see e.g. the textbook \cite{kock2004frobenius} for details.

\subsection{Generalised Verlinde algebras}

To motivate our construction let us consider the simplest case $k=1$ first. Let $z$ be an (invertible) indeterminate and consider the ring $\mathcal{V}_1(n)=\mathbb{C}[z^{\pm 1}][x]/\langle x^n-z\rangle$. The latter has basis $\{1,x,\ldots,x^{n-1}\}$ together with the simple multiplication rule
\[
x^ax^b=\sum_{c=0}^{n-1} z^{\frac{a+b-c}{n}}N_{ab}^cx^c,\qquad N_{ab}^c=\delta_{a+b\text{ mod } n,c}
\]
We have deliberately written the product rule as a sum to facilitate the comparison with the case $k>1$ below. We shall refer to the structure constants $N_{ab}^c$ as {\em fusion coefficients}. Consider the linear map $\varepsilon:\mathcal{V}_1(n)
\to\mathbb{C}[z^{\pm 1}]$ fixed by $\varepsilon(x^a)=\delta_{a,0}/n$. Then for any $z_0\in\mathbb{C}^{\times}$ the quotient $\mathcal{V}_1(n)/(z-z_0)\mathcal{V}_1(n)$ is a Frobenius algebra with the trace map induced by $\varepsilon$. In particular, setting $z_0=1$ the resulting Frobenius algebra is the $\mathfrak{\widehat{sl}}_n$-Verlinde algebra at level $k=1$ \cite{verlinde1988fusion}. It is closely related to the small quantum cohomology of projective space $QH^*(\mathbb{P}^{n-1})$ but both differ as Frobenius algebras as they are endowed with different trace maps.

Based on the work \cite{korffstroppel2010}, which realises the Verlinde algebra for level $k\geq 1$ in terms of a quantum integrable model, $q$-deformed Verlinde algebras were introduced in \cite{korff2013cylindric}. The latter are Frobenius algebras that for $q=0$ specialise to the known $\mathfrak{\widehat{sl}}_n$-Verlinde algebra at level $k$ and for generic $q$ are conjectured \cite{korff2013cylindric}  to be related to Teleman's work \cite{teleman2004k} and \cite{teleman2009index}. Here we consider their $q=1$ specialisation: let $p_r=m_{(r)}$ denote the $r$th power sum and consider the following quotient of $\Lambda_k$,
\begin{equation}\label{V}
\mathcal{V}_k(n)=\mathbb{C}[z^{\pm 1}][x_1,\ldots,x_k]^{S_k}/\langle p_n-zk,p_{n+1}-zp_1,\ldots,p_{n+k-1}-z p_{k-1}\rangle\,.
\end{equation}
For $k=1$ this ring specialises to the above example $\mathcal{V}_1(n)$ and we therefore refer to  \eqref{V} as generalised Verlinde algebra.  

\begin{thm}
(i) The monomial symmetric functions $\{m_\lambda\}_{\lambda\in\mathcal{A}_{k,n}}$ are a basis of  $\mathcal{V}_k(n)$ and together with the linear map $\varepsilon:\mathcal{V}_k(n)\to\mathbb{C}[z^{\pm 1}]$ given by $\varepsilon(m_\lambda)=z^k\delta_{\lambda,n^k}/n^k$ the pair $(\mathcal{V}_k(n),\varepsilon)$ forms a Frobenius extension. (ii) The following product rules hold in $\mathcal{V}_k(n)$,
\begin{equation}\label{fusion}
m_\lambda m_\mu=\sum_{\nu\in\mathcal{A}_{k,n}} z^{\frac{|\lambda|+|\mu|-|\nu|}{n}}N_{\lambda\mu}^{\nu}m_{\nu}\qquad\text{and}\qquad
h_\lambda m_\mu=\sum_{\nu\in\mathcal{A}_{k,n}}z^{\frac{|\lambda|+|\mu|-|\nu|}{n}}\chi_{\nu/d/\mu}(\lambda) m_\nu\;,
\end{equation}
where the fusion coefficients $N_{\lambda\mu}^{\nu}$ coincide with the non-negative integers in \eqref{cylh2h} and $\chi_{\nu/d/\mu}(\lambda)$ is the weighted sum over all cylindric RPP of shape $\nu/d/\mu$ and weight $\lambda$,  i.e. the coefficient of $m_\lambda$ when expanding $h_{\nu/d/\mu}$ defined in \eqref{cylh} into monomial symmetric functions.
\end{thm}
We wish to emphasise that the basis of monomial symmetric functions is distinguished in \eqref{V} by the fact that the fusion coefficients in \eqref{fusion} are non-negative. Other choices, for example Schur functions, lead to structure constants which can be negative.
 
\subsection{Idempotents and the Verlinde formula}
Provided that the roots $t^{\pm1}=z^{\pm 1/n}$ exist, one can show that $\mathcal{V}_k(n)$ is semi-simple over $\mathbb{C}[t^{\pm 1}]$. Assume first that $z=1$ and fix a primitive $n$th root of unity $\zeta\in\mathbb{S}^1\subset\mathbb{C}$. Then each $\alpha\in\mathcal{A}_{k,n}$ determines a point $\zeta^\alpha=(\zeta^{\alpha_1},\ldots,\zeta^{\alpha_k})\in\mathbb{T}^k$ on the $k$-dimensional torus and the polynomials
\begin{equation}\label{idemp}
\mathfrak{e}_{\alpha}(x_1,\ldots,x_k)
=\frac{1}{n^k}\sum_{\lambda\in\mathcal{A}_{k,n}}\frac{|S_\lambda|}{|S_\alpha|}\,
m_{\lambda}(\zeta^{-\alpha})m_\lambda(x_1,\ldots,x_k)
\end{equation}
obey $\mathfrak{e}_{\alpha}(\zeta^\beta)=\delta_{\alpha\beta}$ and form a complete set of idempotents of $\mathcal{V}_k(n)/(z-1)\mathcal{V}_k(n)$. The idempotents for general $z$ are obtained by rescaling the roots, $\zeta^{\pm\alpha_i}\mapsto t^{\pm 1}\zeta^{\pm\alpha_i}$ . 

\begin{thm}
For $z=1$ the product rules \eqref{fusion} describe the multiplication of symmetric functions evaluated at roots of unity and we have the following residue formula for the fusion coefficients in \eqref{cylh2h} and \eqref{fusion}:
\begin{equation}\label{V2}
N_{\lambda\mu}^\nu=\sqrt{\frac{|S_\nu|}{|S_{\lambda}||S_{\mu}|}}\sum_{\alpha\in\mathcal{A}_{k,n}}\frac{\mathcal{S}_{\lambda\alpha}\mathcal{S}_{\mu\alpha}\mathcal{S}_{\nu\alpha}^{-1}}{\mathcal{S}_{\emptyset\alpha}}\;,\qquad \mathcal{S}_{\lambda\alpha}=
\sqrt{\frac{|S_\lambda|}{|S_\alpha|}}
\,\frac{m_{\lambda}(\zeta^\alpha)}{n^{\frac{k}{2}}}\;.
\end{equation}
\end{thm}
Up to the normalisation factor in front of the sum, the last expression \eqref{V2} is a generalistion of the celebrated Verlinde formula \cite{verlinde1988fusion} for $k=1$. The $\mathcal{S}$-matrix in \eqref{V2} links the 2D TQFT $\mathcal{V}_k(n)/(z-1)\mathcal{V}_k(n)$ with the following representation of the modular group:
\begin{prop}
For $z=1$ the matrices $\mathcal{S}_{\lambda\mu}$ and 
$\mathcal{T}_{\lambda\mu}=\delta_{\lambda\mu}\zeta^{-\frac{kn(n-1)}{24}}\theta_{\lambda}$, with $\theta_\lambda=\prod_{i=1}^k\zeta^{\frac{\lambda_i(n-\lambda_i)}{2}}$ and $\lambda,\mu\in\mathcal{A}_{k,n}$, yield a representation of the (double cover of the) modular group $PSL_2(\mathbb{Z})$. That is, we have the identities
\begin{equation}
(\mathcal{ST})^3=\mathcal{S}^2=\mathcal{C},\qquad \mathcal{C}_{\lambda\mu}=\delta_{\lambda\mu^*},
\end{equation}
where $\lambda^*\in\mathcal{A}_{k,n}$ is the image of $(n-\lambda_k,\ldots,n-\lambda_2,n-\lambda_1)$ under $\tau^{m_n(\lambda)}$ and $\mathcal{C}^2=\operatorname{Id}$. In addition, we have unitarity, $\mathcal{S}^{-1}=\mathcal{S}^*$ and the relations  
$\mathcal{S}^*=\mathcal{CS}=\mathcal{SC}$ and $\mathcal{CTC}=\mathcal{T}$. 
\end{prop}
This representation of the modular group is different from the one for the $\mathfrak{\widehat{sl}}_n$-Verlinde algebra in terms of Schur functions if $k>1$; see \cite[Prop 6.15]{kavc1984infinite} and the discussion in \cite[Sec. 6.4]{korffstroppel2010}. To facilitate the comparison, observe that each $\lambda\in\mathcal{A}_{k,n}$ defines a $\mathfrak{\widehat{sl}}_n$-weight $m_n(\lambda)\omega_0+\sum_{i=1}^{n-1}m_i(\lambda)\omega_i$ at level $k$, where the $\omega_i$ are the fundamental $\mathfrak{\widehat{sl}}_n$-weights.

\subsection{The fusion product in terms of the generalised symmetric group}
%
Recall the definition of the generalised symmetric group $S(n,k)=\mathcal{C}_n^{\times k}\rtimes S_k$ as the wreath product of the cyclic group of order $n$ with $S_k$. That is, the group of all complex $k\times k$ matrices that have at most one nonzero entry $\in\mathcal{C}_n$ in each row and column. The simple modules $\mathcal{L}(\boldsymbol{\lambda})$ of $S(n,k)$ are labelled in terms of $n$-multipartitions $\boldsymbol{\lambda}=(\lambda^{(1)},\ldots,\lambda^{(n)})$ of $k$, i.e. $|\boldsymbol{\lambda}|=\sum_i|\lambda^{(i)}|=k$; see e.g. \cite{osima1954representations} for details. Consider the following exact sequence of groups, %
$
1\rightarrow \mathcal{C}_{n}^{\times k}\hookrightarrow
S(n,k)\twoheadrightarrow S_{k}\rightarrow 1,$ %
and denote by $y_i$ the image of the generator of the $i$th copy of $\mathcal{C}_n$ of the normal subgroup $\mathcal{N}\cong\mathcal{C}_{n}^{\times k}$ in $S(n,k)$.
 \begin{lem}[\cite{osima1954representations}]
The characters of $\mathcal{L}(\boldsymbol{\lambda})$ restricted to the normal subgroup $\mathcal{N}$ are given by
\begin{equation}
\chi_{\boldsymbol{\lambda}} (y^{\alpha}) =\operatorname{Tr}_{\mathcal{L}(\boldsymbol{\lambda})}y^\alpha=f_{\boldsymbol{\lambda}}\, m_{\lambda}(\zeta^{\alpha}),\quad\forall\alpha\in\mathbb{Z}_n^{\times k},
\quad f_{\boldsymbol{\lambda}}=\prod_{i=1}^{n} f_{\lambda^{(i)}}\,, 
\end{equation}
where $f_{\lambda^{(i)}}$ is the number of standard tableaux of shape $\lambda^{(i)}$ and $\lambda\in\mathcal{A}_{k,n}$ is the unique partition with $m_i(\lambda)=|\lambda^{(i)}|$.
\end{lem}

 Let $\operatorname{Rep}S(n,k)$ denote the representation ring and $c_{\boldsymbol{\lambda}\,\boldsymbol{\mu}}^{\boldsymbol{\nu}}$ its structure constants with respect to the simple modules $\mathcal{L}(\boldsymbol{\lambda})$. 
 
 \begin{cor} 
 We have the following alternative expression for the fusion coefficients in \eqref{cylh2h}, \eqref{reducedN} and \eqref{fusion},
\begin{equation}
N_{\lambda\mu}^{\nu}=\sum_{\boldsymbol{\nu}}c_{\boldsymbol{\lambda}\,\boldsymbol{\mu}}^{\boldsymbol{\nu}}
\,\frac{f_{\boldsymbol{\nu}}}{f_{\boldsymbol{\lambda}}f_{\boldsymbol{\mu}}},
\end{equation} 
where the sum runs over all multipartitions $\boldsymbol{\nu}$ such that $|\nu^{(i)}|= m_i(\nu)$ with $\nu\in\mathcal{A}_{k,n}$. Similarly, $\boldsymbol{\lambda}$, $\boldsymbol{\mu}$ is any pair of multipartitions which satisfy the analogous condition for $\lambda,\mu\in\mathcal{A}_{k,n}$.
\end{cor}
We can interpret this result as stating that the generalised Verlinde algebra \eqref{V} describes for $z=1$ the representation ring of the normal subgroup $\mathcal{N}$ inside $\operatorname{Rep}S(n,k)$.

\subsection{Schur-Weyl duality and the principal subalgebra}
Define the following elements of the loop algebra $\mathfrak{sl}_n[z,z^{-1}]= \mathfrak{sl}_{n}(\mathbb{C})\otimes\mathbb{C}[z,z^{-1}]$ in terms of the unit matrices $\{e_{ij}~|~1\leq i,j\leq n\}\subset\mathfrak{gl}_{n}(\mathbb{C})$,
\begin{equation}\label{Pr}
P _{r}=\sum_{i-j=r}e_{ij}+z\sum_{j-i=n-r}e_{ij}\;,\qquad r=1,\ldots ,n-1,
\end{equation}%
and set $ P_{r+n}=z P_{r}$ for all other $r\in \mathbb{Z}\backslash n\mathbb{Z}$. The latter generate a commutative subalgebra of $\mathfrak{sl}_n[z,z^{-1}]$ with respect to the usual Lie bracket, $[f(z)x,g(z)y]:=f(z)g(z)[x,y]$, where $f,g\in\mathbb{C}[z,z^{-1}]$ and $x,y\in\mathfrak{sl}_n(\mathbb{C})$. This commutative subalgebra is the image of the {\em principal Heisenberg subalgebra} in the affine Lie algebra $\mathfrak{\widehat{gl}}_n=\mathfrak{gl}_n[z,z^{-1}]\oplus\mathbb{C}\underline{k}$ under the projection defined via  $0\rightarrow\mathbb{C}\underline{k}\rightarrow \mathfrak{\widehat{gl}}_n \overset{\pi}{\rightarrow} \mathfrak{gl}_n[z,z^{-1}]\rightarrow 0$, where $\underline{k}$ is the central element; see \cite{kac1994infinite} for further details.

Recall the canonical left action of the universal enveloping algebra $U=U(\mathfrak{sl}_n[z,z^{-1}])$ on the tensor product $\mathbb{C}[z,z^{-1}]\otimes (\mathbb{C}^n)^{\otimes k}$: let $x\in \mathfrak{sl}_n$ and $v_1,\ldots,v_k\in\mathbb{C}^n$, then
\[
 f(z)\otimes x.g(z)v_1\otimes\cdots \otimes v_k=f(z)g(z)\sum_{i=1}^k v_1\otimes\cdots\otimes x.v_i\otimes \cdots\otimes v_k
\;.
\] 
This $U$-action commutes with the canonical right action of the symmetric group $S_k$ on $(\mathbb{C}^n)^{\otimes k}$. Thus, the $U$-action leaves the $k$th symmetric power $\mathbb{C}[z,z^{-1}]\otimes S^k(\mathbb{C}^n)$ invariant.
\begin{thm} 
Denote by $\mathbb{V}_k(n)\subset\operatorname{End}(\mathbb{C}[z^{\pm 1}]\otimes S^k(\mathbb{C}^n))$ the commutative ring generated by the $\{P_r\}_{r\in\mathbb{Z}}$ with $P_{mk}=z^m k\,\operatorname{Id}$. Then the map $\mathbb{V}_k(n)\to\mathcal{V}_k(n)$ fixed by $P_r\mapsto p_r$ and $P_{-r}\mapsto z^{-1}p_{n-r}$ for $r=1,\ldots,n-1$, where $p_r=m_{(r)}$ are the power sums, is a ring isomorphism. 
\end{thm}
The analogous statement for the antisymmetric power $\bigwedge^k(\mathbb{C}^n)$ yields a ring isomorphism with the small quantum cohomology ring $QH^*(\operatorname{Gr}_k(\mathbb{C}^n))$ of the Grassmannian, where $z=(-1)^{k-1}q$ and $q$ is the deformation parameter \cite{KorffPalazzo}. The last theorem shows that the generalised Verlinde algebra \eqref{V} naturally arises as a counterpart of $QH^*(\operatorname{Gr}_k(\mathbb{C}^n))$ within the context of Schur-Weyl duality and that in this setting the cylindric complete symmetric functions \eqref{cylh} play a role analogous to Postnikov's toric Schur functions \cite{postnikov2005affine}.


\end{document}